\documentclass[a4paper,12pt]{article}
\usepackage{a4wide}
\usepackage{graphicx,amsbsy,color}
\usepackage{amsbsy,amsmath, amsfonts}
\usepackage{amsthm,amssymb}
\usepackage{verbatim}   % useful for program listings
\usepackage{color,xcolor,tikz,pgfplots}% use if color is used in text
\pgfplotsset{compat=1.7}%\usepackage{subfigure}  % use for side-by-side figures
\usepackage{amssymb}
\usepackage{hyperref}
\usepackage{caption}
\usepackage{subcaption}
\usepackage{fancyhdr}
\usepackage{lineno,hyperref}
\usepackage{float}
\usepackage{url}
\usepackage[sort&compress]{natbib}
\DeclareMathOperator*{\argmin}{argmin}

\theoremstyle{definition}
\newtheorem{rem}{Remark}
\usepackage{numcompress}
\usepackage{etoolbox}

\makeatletter

\pretocmd{\NAT@citex}{%
  \let\NAT@hyper@\NAT@hyper@citex
  \def\NAT@postnote{#2}%
  \setcounter{NAT@total@cites}{0}%
  \setcounter{NAT@count@cites}{0}%
  \forcsvlist{\stepcounter{NAT@total@cites}\@gobble}{#3}}{}{}
\newcounter{NAT@total@cites}
\newcounter{NAT@count@cites}
\def\NAT@postnote{}

\def\NAT@hyper@citex#1{%
  \stepcounter{NAT@count@cites}%
  \hyper@natlinkstart{\@citeb\@extra@b@citeb}#1%
  \ifnumequal{\value{NAT@count@cites}}{\value{NAT@total@cites}}
    {\ifNAT@swa\else\if*\NAT@postnote*\else%
     \NAT@cmt\NAT@postnote\global\def\NAT@postnote{}\fi\fi}{}%
  \ifNAT@swa\else\if\relax\NAT@date\relax
  \else\NAT@@close\global\let\NAT@nm\@empty\fi\fi% avoid compact citations
  \hyper@natlinkend}
\renewcommand\hyper@natlinkbreak[2]{#1}

\patchcmd{\NAT@citex}
  {\ifNAT@swa\else\if*#2*\else\NAT@cmt#2\fi
   \if\relax\NAT@date\relax\else\NAT@@close\fi\fi}{}{}{}
\patchcmd{\NAT@citex}
  {\if\relax\NAT@date\relax\NAT@def@citea\else\NAT@def@citea@close\fi}
  {\if\relax\NAT@date\relax\NAT@def@citea\else\NAT@def@citea@space\fi}{}{}

\makeatother
\graphicspath{
{./}
{./Pics/}
{./Pics/MeshAdapt/}
{./Pics/MeshAdapt/MatLab/}
{./PIcs/}
{./pIcs/MeshAdapt/}
{./pics/}
}

\begin{document}

%\input title_page %First page
%\clearemptydoublepage

\pagenumbering{arabic}

\title{\vspace*{-4.em} Adaptive Eigenspace Segmentation}
\author{Uri Nahum and Philippe C. Cattin\\ Department of Biomedical Engineering, \\ 
University of Basel, Allschwil, Switzerland}
\date{\today}

\maketitle
\begin{abstract}
Image segmentation is an inherently ill-posed problem and thus requires regularization in order to limit the search space to \emph{reasonable} solutions. A majority of segmentation methods integrates these regularization terms in one way or the other in an energy functional using a balancing term. The tuning of this parameter that either favours more the regularization or the data conformity is critical and, unfortunately, the success of the optimization process strongly depends on it. Often the optimal settings change from image to image.

In this paper we propose a novel general framework based on an adaptive eigenspace that was first proposed for solving inverse problems. 
The resulting method proves accurate and yields robust results, without the need for optimization techniques or being 
sensitive to the parameter choice. In fact, the method solves a symmetric positive definite sparse system and hence,
uses only a fraction of the computational cost. The method is very versatile and does not need parameter-tuning,
when segmenting objects from any kind of an image or when segmenting different organs. As the adaptive eigenspace is determined directly from the image to segment, the approach also does not need a tedious training phase.
\end{abstract}

%\tableofcontents

%\newpage
\section{Introduction}
\label{s:intro}

Segmentation  is a basic image processing technique that has spurred widespread interest during the last decades. Many of the proposed approaches use some form of regularization and often employ an iterative scheme. Regularization is recommended as it limits the result space to “meaningful” solutions.
 
Image segmentation is the process in which an algorithm divides a digital image into groups of connected pixels.
The idea is to split the image domain into non-intersecting segments and that can subsequently be analyzed to extract higher-level information from the image.
Segmentation of images has been an area of active research over more than 45\,years now \cite{pavlidis1972segmentation} and yet, finding a robust algorithm,
which is able to segment different types of images without intensive parameter tuning or training is still a challenge. Furthermore, measurement noise and speckle, as seen in ultrasound or Optical Coherence Tomography (OCT) images, still pose a largely unsolved challenge.
Many techniques were proposed over the years, see \cite{Angenent2006} for a concise overview, including very basic ones, such as thresholding, region growing or watershed, see for example \cite{adams1994seeded}, \cite{fan2001automatic} and \cite{shafarenko1997automatic}.
The idea behind region growing, for example, is to start in some seed points and to test if the neighboring pixels should be part of their segment.
Those methods are very intuitive, but are not robust under noise.

In recent years, segmentation using deep learning and neural networks has become very popular and many papers were written and are still written on this topic,
for example \cite{schmidhuber2015deep}, \cite{litjens2017survey} and \cite{andermatt2016multi}. 
%Although, the fundamental theory behind these approaches has been known for decades \cite{LeCun1998} only the recent developments in GPU hardware \cite{Steinkraus2005} made the computationally complex training phase tractable. %Given enough labelled training data is available, the deep learning approaches deliver show a big versatility and show often very good results close to human readers.
In contrast to the above mentioned methods, where the information is taken from the image to be segmented, in deep learning, 
the method tries to automatically learn the segmentation through a large dataset-collection of sample images and their segmentation's correct labelling.% (often done manually).
%In other words, the method finds the optimal weights for the network through the large number of examples. 
%After the initial training data phase, the method applies the weights to new unseen images for segmentation.

Besides the basic ad hoc segmentation approaches described above, a family of more principled segmentation approaches form the energy based techniques. This family includes techniques such as {\it active contours}, also known as {\it snakes}, see \cite{kass1988snakes}
and its various extensions \cite{cohen1991active,Kichenassamy1996,caselles1997geodesic,xu19972d,bresson2007fast}.	
In this method, a snake is a spline, which can be deformed to minimize the energy functional.
The snake is influenced by constraints, called external and internal forces that deform it to match the contours of an object. Some of these energies force the snake to the contour of the image while others regularize the resulting contour to a “reasonable” result. Besides the snakes, other energy based techniques became popular such as the level-sets \cite{Chan2001} or the graph-cut \cite{BOYKOV2001} approaches. All these techniques represent different ways of solving an energy term which is similar for all of them. In fact they all try to minimize a positive cost functional $E$, where $I$ is the image to segment and $u$ the segmentation result
\begin{equation}\label{basicOpt}
 E[u;I] = F[u;I] + \eta \mathcal{R}[u]\,.
 \end{equation}

The energy $E$ is composed of two components, namely the \emph{fidelity term} $F$ that forces the segmentation
to match the source image as close as possible and the \emph{regularization term} $\mathcal{R}$ ensures that the resulting segmentation is reasonable.
The balancing term $\eta$ is common to these methods and is a parameter that needs careful tuning to each and every application separately as it balances the trade-off between closeness of the solution to the image and the regularization term.
Depending on the application one may choose different kinds of penalty functionals $\mathcal{R}$ for regularization, including the common $L^2$-norm, the {\it penalized TV-regularization},
the {\it Sobolev $H^1$-penalty} functional, the {\it Lorentzian} penalty functional~\cite{Tikhonov:1943:SIP,Gene:1999:TRT,Engl:1996:RIP}.

The {\it adaptive eigenspace} (AE), initially introduced in \cite{deBuhan:2009:SRP} and later in \cite{deBuhan:2013:AIM} and in \cite{de2017numerical} which is strongly related to \cite{gilboa2016nonlinear}, was used in~\cite{Grote:2016:AEI} and in \cite{grote2017adaptive} as a new regularization method,
where the AE of the penalized TV-gradient was used to regularize an inverse medium problem.
 Instead of adding a penalty functional $\mathcal{R}$ to the minimization problem, we build the parameter from the eigenfunctions of the 
 TV-regularization gradient drastically simplifying the optimization of the energy functional. The method was able to reduce dramatically the number of variables and to achieve very accurate solutions 
 in smaller computational cost.
 In \cite{Nahum:2016:AEI}, several eigenspaces from different regularization terms are introduced, one of those was the adaptive eigenspace of the Lorentzian regularization.
 In \cite{gilboa2016nonlinear}, the nonlinear spectral representation is introduced.
 There, eigenfunctions of the TV-regularization and other convex functionals are used for image decomposition and denoising.
 The idea of segmenting images using eigenfunctions is also introduced in \cite{7433409},
 where eigenfunctions of an anisotropic diffusion operator were used successfully for image segmentation.
 
In this paper, we propose a \emph{general and versatile framework} using the adaptive eigenspace of the non-convex Lorentzian penalty functional for segmentation. 
Instead of optimizing over the energy functional $E$,
we compute the eigenfunctions of the gradient of the regularization term $\mathcal{R}$ and find the segmentation there.
Here, we show how this approach can segment useful information with very low computational cost.
As shown in~\cite{Nahum:2016:AEI}, the method does not need any parameter tuning, when recovering different media, the same applies for image segmentation.
Using the AE for segmentation yields a fast, sparse and reliable method, which is very robust under noise or speckle. As the approach has no parameters there is no need for parameter-tuning such as balancing the term $\eta$ from \eqref{basicOpt} and is thus inherently insensitive to parameter tuning. Lastly, we show that the proposed novel segmentation approach can also be applied for noise reduction by \textit{segmenting the noise out} (rather than filtering the noise).

This paper is organized in four parts. In Section \ref{s:Segmentation}, we show how to derive the AE of the Lorentzian.
In Section \ref{s:AEI_Principle} properties of the AE are introduced.
In Section \ref{s:PpT}, we show mathematical evidence of the robustness and usefulness of the AE, using examples in one or two space dimensions. 
At last, in Section \ref{s:medical} we show some numerical results that underpin the efficiency and suitability of the AE to segment images, remove noise and combinations thereof.

\section{Image Segmentation and Regularization}
\label{s:Segmentation}

% Image processing and segmentation is a wide area, which uses mostly iterative methods and regularization techniques
% to remove noise or to separate elements out of an image.
% More than 25 years now methods like penalized {\it Total Variation} (TV) regularization~\cite{Rudin:1992:NTV},
% are used for noise removal in digital image processing. The TV regularization can reconstruct media with discontinuities; 
% it is able to remove unwanted noise, while preserving important details.
% In inverse problems, adding a penalty-terms such as TV-regularization are called
% Tikhonov regularization \cite{Tikhonov:1943:SIP,Gene:1999:TRT,Engl:1996:RIP} and can tackle ill-posedness\cite{maharramov:2014:RFW,Esser:2015:TVR}.
% Using the TV-regularization can reconstruct nearly piecewise constant images with high quality \cite{Dobson:1996:RBI}. 

% In ~\cite{Grote:2016:AEI}, a new method to regularize the inverse medium problem was introduced.
% There, the adaptive eigenspace (AE) of the penalized TV-gradient was used to regularize the inverse medium problem.
% Instead of adding penalty functional to the minimization problem, we build the parameter from the eigenfunctions of the 
% TV-regularization gradient. The method was able to reduce dramatically the number of variables and to achieve very accurate solutions 
% in smaller computational cost.
% In \cite{Nahum:2016:AEI} several eigenspaces from several regularization terms are introduced, 
% one of those was the adaptive eigenspace of the Lorenzian regularization.

In image segmentation one uses, usually, a minimization problem of a positive energy functional $E$ to extract wanted information out of an image. For a given image $I$
and set of admissible parameters $U_{ad}$, one solves
\begin{equation}\label{optimization}
 u^*=\argmin\limits_{u\in U_{ad}} \left\{ F(u,I)+\eta\mathcal{R}(u)\right\}\,,
\end{equation}
where $F$ is the fidelity term and $\mathcal{R}$ is a regularization functional, which provides any additional knowledge on the wished segmentation
and tackles ill-posedness of the optimization problem. The parameter $\eta$ balances between those functionals.

For the regularization, one can choose different penalty functionals $\mathcal{R}(u)$:
the standard Tikhonov~\cite{Tikhonov:1943:SIP,Vogel:2002:CMIP,Gene:1999:TRT,Engl:1996:RIP} $L^2${\it-penalty} functional 
\begin{equation}\label{eq:L2reg}
 \mathcal{R}_{L^2}(u)=\frac{1}{2}\|u\|^2_{L^2}\,,
\end{equation}
or the {\it Sobolev $H^1$-penalty} functional
\begin{equation}\label{eq:H1reg}
 \mathcal{R}_{\nabla u}(u)=\frac{1}{2}\int_\Omega \sum_{i=1}^d \left(\dfrac{\partial u}{\partial x_i}\right)^2\,dx\,,
\end{equation}
which penalizes strong variation in the solution and leads to smooth results. 
% These two penalty functionals, can be generalized \cite{Gene:1999:TRT} by 
% \begin{equation}\label{eq:generalTikForm}
%  \mathcal{R}_D(u)=\frac12\|D\,u \|^2_{L^2} \,. 
% \end{equation}
% Usually, $D$ is a first or second space derivative operator, but other choices are possible.

A popular penalty functional penalizes the Total Variation (TV) \cite{Vogel:1996:IMTV},
and was originally introduced for noise removal in image processing \cite{Rudin:1992:NTV}.
While preserving important detail, regularizing images using the penalized TV, removed unwanted noise.
It uses the $L^1$ norm and is given by
\begin{equation}\label{eq:TVreg}
 R_{TV}(u) = \frac12\int_\Omega \sqrt{|\nabla u|^2+\varepsilon^2}\,dx\,,\quad \varepsilon\neq0\,.
\end{equation}
% Note that the penalized TV-regularization functional is differentiable. 

An important penalty functional from a group of non-convex regularization terms is the Lorentzian penalty term 
\begin{equation}
R_{Lorentz}(u)= \frac12\int_{\Omega} \dfrac{\gamma|\nabla u|^2}{1+\gamma|\nabla u|^2}\,dx\,,\quad \gamma>0\,.
\end{equation}
It penalizes
strong variations in the solution and contains an extra parameter $\gamma$ to allow discontinuities.
%For simplicity we use the notation $R_{TV}$ to denote
%the penalized TV functional and $\mathcal{R}_{\nabla u}$ for the non-penalized TV functional \eqref{eq:H1reg}.
% This regularization penalizes the gradient of $u$ and allows discontinuities in the profile through
% the penalization parameter $\varepsilon$ as well.

In contrast, when segmenting an image with the proposed adaptive eigenspace we do not optimize \eqref{optimization} but only have to compute the image's eigenfunctions of the gradient of the regularization term $\mathcal{R}$.
This approach has proved itself as very accurate and reliable in the inverse medium problems in \cite{Grote:2016:AEI}. There, a severely ill-posed problem was solved and regularized using the adaptive eigenspace of the TV-regularization
gradient with much success.

Unlike the inverse medium problem, here, the image to be segmented is known and we may use the image itself to build the adaptive eigenspace. 
Hence, we can use an AE of non-convex penalty terms, which may be more sensitive to changes in the gradient of the image.
In \cite{Nahum:2016:AEI} and in \cite{grote2017adaptive}, several eigenspaces of different regularization terms are introduced, 
one of those was the adaptive eigenspace of the Lorentzian regularization.
Here, we may use this regularization without intensifying the ill-posedness or non-convexity of the problem.

This results in an image-processing approach, which is not sensitive to a parameter choice and there is no need to change parameters when segmenting different kinds of images, including photos or medical images.

	\section{Principle of the Adaptive Eigenspace (AE)}
	\label{s:AEI_Principle}
  
The main purpose of image segmentation using the AE lies on the parametrization of the image~$I$.
Standard techniques use regularization methods and iterative methods on the image. Following~\cite{Grote:2016:AEI}, 
we propose to use the adaptive eigenspace as regularization.
There, the parameter is unknown and we use the AE of the TV-regularization as regularization.
Here, as the image $I$ is known, we follow~\cite{Nahum:2016:AEI} and derive the adaptive eigenspace of the Lorentzian regularization.
For an image $I$ the Lorentzian regularization-functional is given by 
\begin{equation}\label{eq:Lorentzian_reg}
R_{Lorentz}(I)= \frac12\int_{\Omega} \dfrac{\gamma|\nabla I(x)|^2}{1+\gamma|\nabla I(x)|^2}\,dx\,,\quad \gamma>0\,.
\end{equation}
Next, we compute the gradient of \eqref{eq:Lorentzian_reg}
\begin{equation}\label{eq:LorentzRegGrad}
    \nabla_u \mathcal{R}_{Lorentz}(I)=  -\nabla\cdot
    \left(\dfrac{\gamma\nabla I(x)}{\left(1+\gamma|\nabla I(x)|^2\right)^2}\right)\,,\quad \gamma>0\,.
    \end{equation}    
Following \cite{Nahum:2016:AEI}, we build an adaptive eigenspace from \eqref{eq:LorentzRegGrad}, this means that
we take the gradient from \eqref{eq:LorentzRegGrad}, substitute $I(x)$ into $\phi(x)$ only where $I(x)$ does not appear in an absolute value.
Then, we impose Dirichlet boundary conditions and get the following eigenspace problem
\begin{equation} \label{eq:eigenfunctionsTV}
	\left \lbrace \begin{array}{rclll}
		- \nabla \cdot \left(\dfrac{\gamma\nabla \phi_m(x)}{{\left(1+\gamma|\nabla I(x)|^2\right)^2}} \right) & = & \lambda_m \phi_m(x), & \qquad & \forall\,x\in\Omega, \\[0.2em]
		\phi_m(x) & = & 0, & & \forall\,x\in\Gamma\,.
	\end{array} \right.
\end{equation}
Now, we can expand $I(x)$ as
\begin{equation} \label{eq:param_u}
	I \ = \ I_0(x) \ + \ \sum_{m\geq 1} \beta_m \phi_m(x),
\end{equation}
where $I_0 \in H^1(\Omega)$ is a prolongation of~$c^2_{|\Gamma} \in H^\frac{1}{2}(\Gamma)$ and the 
functions $\phi_m \in H_0^1(\Omega)$ form a Hilbertian basis to parametrize~$I-I_0$. 
The prologation~$I_0$ is the eigenfunction corresponding to $\lambda = 0$, which holds the same boundary data as $I$, namely
\begin{equation} \label{eq:u_0}
	\left \lbrace \begin{array}{rclll}
		- \nabla \cdot \left(\dfrac{\gamma\nabla I_0(x)}{{\left(1+\gamma|\nabla I(x)|^2\right)^2}} \right) & = & 0, & \qquad & \forall\,x\in\Omega, \\[0.2em]
		I_0(x) & = & I(x), & & \forall\,x\in\Gamma.
	\end{array} \right.
\end{equation}
\section{Properties of the Adaptive Eigenspace Basis}
\label{s:PpT}
	\subsection{One-dimensional case}
	\label{ss:PpT_1D}

In \eqref{eq:param_u}, we have used the basis of eigenfunctions $\{\phi_m \}_{m\geq 1}$ defined 
by~\eqref{eq:eigenfunctionsTV} together with $I_0$ defined by \eqref{eq:u_0}. In this section, for a given
image $ I(x) $, we provide some analytical and numerical evidence which underpins the basis choice for segmentation of images with or without noise. Similar examinations for proving the usefulness of the methods were made for the penalized TV-regularization eigenspace in \cite{grote2017adaptive} and for the spectral TV in \cite{gilboa2016nonlinear}.

We can approximate an image $I$ in the eigenspace $\{\phi_1,\phi_2,\,\ldots,\,\phi_K\}$, where 
all~$\phi_m(x)$ satisfy~\eqref{eq:EF-1D} in one space dimension:
\begin{equation} \label{eq:EF-1D}
	\left \lbrace \begin{array}{c}
		- \dfrac{d}{dx} \left( \mu(x) \dfrac{d}{dx} \phi_m(x) \right)  =  \lambda_m\phi_m(x)  \qquad  \forall\,x\in [a,b], \\[0.2em]
		\phi_m(a)=0\,,\quad \phi_m(b)=0,
	\end{array} \right.
\end{equation}
where
\begin{equation} \label{eq:mu-1D}
	\mu(x)=\dfrac{\gamma}{(|1+\gamma I^\prime(x)^2)^2}\,, \quad \forall x\in [a,b].
\end{equation}
The behavior of $\phi_m(x)$ strongly depends on the magnitude of $I^\prime(x)$.

If $|I^\prime(x)|\simeq C > \frac{1}{\sqrt{\gamma}}$ in some part of~$\Omega,$ $\phi_m(x)$ behaves like 
$\phi_m(x)\simeq A \sin(C^2\sqrt{\lambda_m}\,x)+ B \cos(C^2\sqrt{\lambda_m}\,x)$. However, if~$I$ is  
constant in some part of~$\Omega$, $|I^\prime(x)|= 0$, then $\mu=\gamma$ there, and $\phi_m$  
behaves like $\phi_m(x)\simeq A_m \sin(\sqrt{\frac1\gamma\lambda_m}\,x)+ B_m \cos(\sqrt{\frac1\gamma\lambda_m}\,x)$.
For large enough $\gamma$, $\phi_m(x)$ has very slow variation and remains constant.

\begin{figure}[p]
	\centering
	\vspace*{-1.truecm}
	\includegraphics[width=0.35\textwidth]{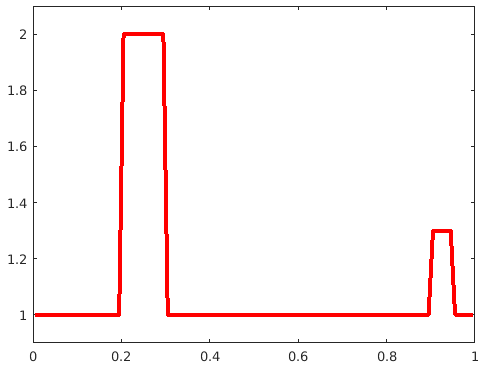}\\
	\includegraphics[width=0.35\textwidth]{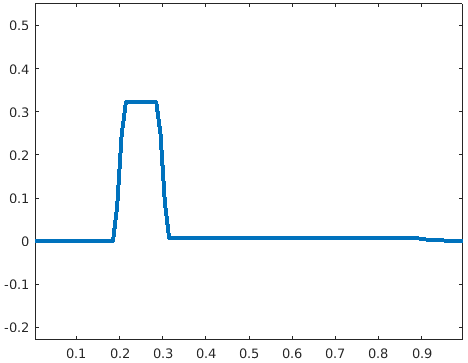}\hspace*{1.truecm}
	\includegraphics[width=0.35\textwidth]{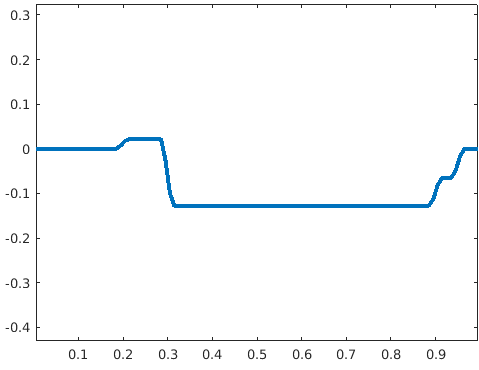}\\
	\includegraphics[width=0.35\textwidth]{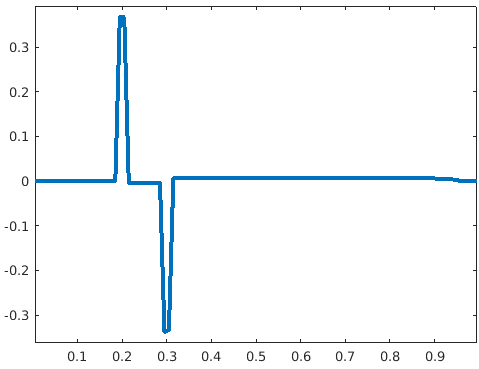}\hspace*{1.truecm}
	\includegraphics[width=0.35\textwidth]{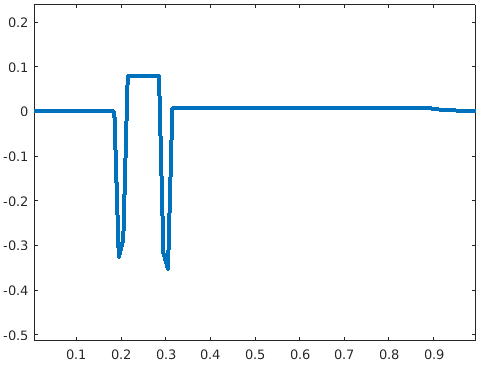}\\
	\includegraphics[width=0.35\textwidth]{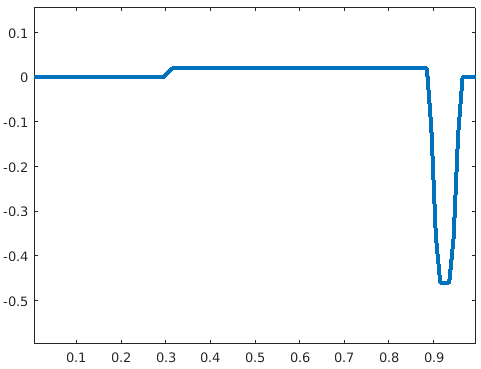}\hspace*{1.truecm}
	\includegraphics[width=0.35\textwidth]{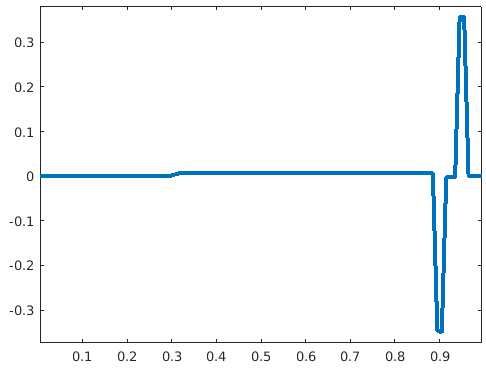}\\
	\includegraphics[width=0.35\textwidth]{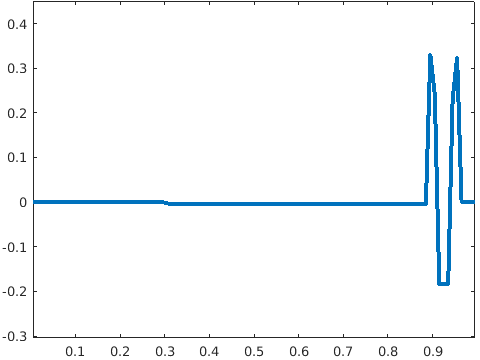}\hspace*{1.truecm}
	\includegraphics[width=0.35\textwidth]{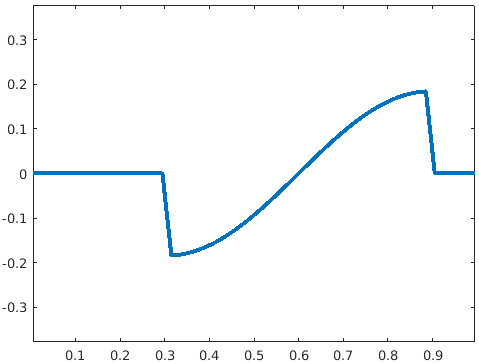}
	\caption{One-dimensional case. The true profile $I$ (top), together with the first eight eigenfunctions from~\eqref{eq:EF-1D} with~$\gamma = \max(|I^\prime|)$. \label{fig:phi-1D}}
\end{figure}

To illustrate this behavior, we now consider the profile~$I(x)$ containing two one-dimensional objects, as shown on the top of Fig.~\ref{fig:phi-1D}, note that 
$I^\prime(x) \neq 0$ for $0.19 \leq x \leq 0.2$, $0.29 \leq x \leq 0.3$, $0.89 \leq x \leq 0.9$ and $0.95 \leq x \leq 0.96$ and $I^\prime(x) = 0$, otherwise. In Fig.~\ref{fig:phi-1D},
we present $I(x)$ together with some of its 
eigenfunctions from problem~\eqref{eq:EF-1D}. We can see that $\phi_1(x)$ coincides to the first object in $I(x)$ and $\phi_5(x)$ to the second object in $I(x)$ (up to a constant).
Actually, all the eigenfunctions, including those seen in Fig. \ref{fig:phi-1D}, can be described as 
\begin{equation}
  \phi_m(x)\simeq  \begin{cases}
    A_{m,1} \sin(\frac{1}{\sqrt\gamma}\sqrt{\lambda_m}\,x)+ B_{m,1} \cos(\frac{1}{\sqrt\gamma}\sqrt{\lambda_m}\,x), & x\in [0.0,0.19), \\
    A_{m,2} \sin(10^4\sqrt{\lambda_m}\,x)+ B_{m,2} \cos(10^4\sqrt{\lambda_m}\,x), & x\in [0.19,0.2), \\
    A_{m,3} \sin(\frac{1}{\sqrt\gamma}\sqrt{\lambda_m}\,x)+ B_{m,3} \cos(\frac{1}{\sqrt\gamma}\sqrt{\lambda_m}\,x), & x\in [0.2,0.3), \\
    A_{m,4} \sin(10^4\sqrt{\lambda_m}\,x)+ B_{m,4} \cos(10^4\sqrt{\lambda_m}\,x), & x\in [0.3,0.31), \\
A_{m,5} \sin(\frac{1}{\sqrt\gamma}\sqrt{\lambda_m}\,x)+ B_{m,5} \cos(\frac{1}{\sqrt\gamma}\sqrt{\lambda_m}\,x), & x\in [0.31,0.89)\\
A_{m,6} \sin(30^2\sqrt{\lambda_m}\,x)+ B_{m,6} \cos(30^2\sqrt{\lambda_m}\,x), & x\in [0.89,0.9]\\
A_{m,7} \sin(\frac{1}{\sqrt\gamma}\sqrt{\lambda_m}\,x)+ B_{m,7} \cos(\frac{1}{\sqrt\gamma}\sqrt{\lambda_m}\,x), & x\in [0.9,0.95)\\
A_{m,8} \sin(30^2\sqrt{\lambda_m}\,x)+ B_{m,8} \cos(30^2\sqrt{\lambda_m}\,x), & x\in [0.95,0.96)\\
A_{m,9} \sin(\frac{1}{\sqrt\gamma}\sqrt{\lambda_m}\,x)+ B_{m,9} \cos(\frac{1}{\sqrt\gamma}\sqrt{\lambda_m}\,x), & x\in [0.96,1]. 
 \end{cases}
\end{equation}
For every interval, $\phi_m(x)$ has a different frequency, which is defined by the expression~$C^2\sqrt{\lambda_m}$. 
For example, for each eigenfunction $\phi_m(x)$, the frequency in the subinterval~$[0.3,0.31]$ is $\left(\frac{10}{3}\right)^2$ times higher than the 
frequency in the subinterval~$[0.95,0.96)$ due to the strong dependency of $\mu$ from \eqref{eq:mu-1D} in $I^\prime$ for $I^\prime\neq 0$ (it appears in high potency in the divisor of $\mu$). 
In the subintervals~$[0.0,0.19),[0.2,0.3),[0.31,0.89),[0.9,0.95)$ and $[0.96,1)$, the frequency depends strongly 
on~$\gamma$ which is chosen as $\gamma=\max(|I^\prime(x)|$. While $\frac{1}{\sqrt\gamma}\sqrt{\lambda_m}$ is 
very small, the eigenfunction~$\phi_m(x)$ oscillates very slowly in this subinterval and thus behaves as a constant. 
As $\lambda_m$ increases, the frequency $\frac{1}{\sqrt\gamma}\sqrt{\lambda_m}$ increases as well and oscillations appear (see $\phi_7$ in Fig. \ref{fig:phi-1D}, bottom right). 
Clearly, for high enough $\gamma$, more eigenfunctions $\phi_m(x)$ essentially behave as constants on this 
subinterval.
\bigskip

\begin{figure}[!ht]
	\centering
	\includegraphics[width=0.32\textwidth]{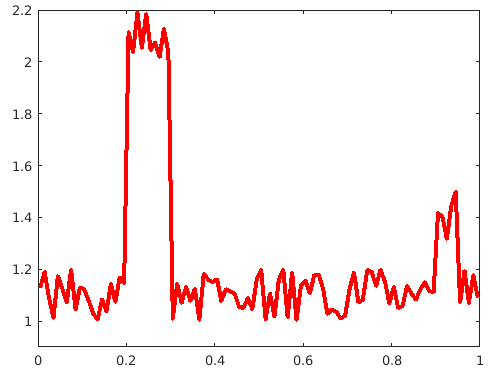}
	\includegraphics[width=0.32\textwidth]{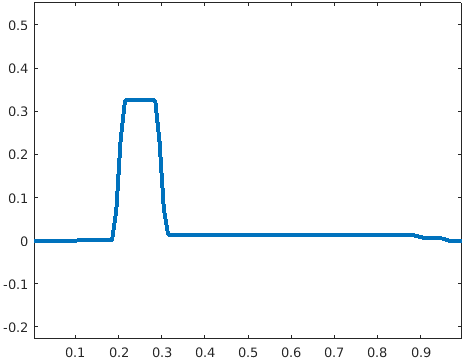}
	\includegraphics[width=0.32\textwidth]{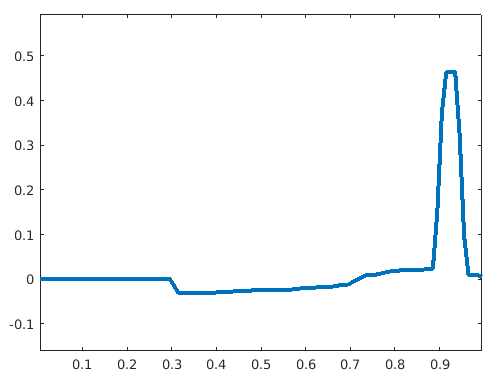}
	\caption{One-dimensional case. The noised image~$I_{noise}$ (left) and 
	the first and fifth eigenfunctions from~\eqref{eq:EF-1D}~$\phi_1$ (center) and $\phi_5$ (right).}
	\label{fig:phi_glatt-1D}
\end{figure}

Finally, we consider~$I$, shown in the top of Fig.~\ref{fig:phi-1D}, with $20\%$ of added noise such that 
$$ I_{noise} = I\left(1+\delta\,\xi\right)\, $$
where $\xi$ is uniformly distributed random number in the interval $(0,1)$ and $\delta=0.2$ represents the noise level.
Since $\mu$ is strongly dependent in $I^\prime$, for $I^\prime\neq 0$, the first eigenfunction will essentially extract the objects from the added noise. 
%This can be done as long as the noise is not to big to dominate the high derivatives on the boundary of the objects in the image.
The image $I_{noise}$ with $20\%$ of added noise is shown in Fig~\ref{fig:phi_glatt-1D}, together with the first eigenfunction~$\phi_1$ and the fifth eigenfunction~$\phi_5$ obtained from~\eqref{eq:EF-1D} using $I_{noise}$. 
Again, we observe that two of the eigenfunctions captures the elements appears in $I_{noise}$, this is up to a small perturbation, as can be seen in both eigenfunctions.
\begin{rem}\label{rem:thresholding}
although the eigenfunctions are not strictly segmentations as the are not binary labels, one can simply threshold them to yield binary segmentations.
%The minor perturbations in the eigenfunctions can be easily removed by using a simple sparsity strategy, for example in \cite{Nahum:2016:AEI,Grote:2016:AEI}.

\end{rem}

	\subsection{Two-dimensional case}
	\label{ss:PpT_Laplace}

\begin{figure}[h!]
	\centering	
	\includegraphics[width=0.32\textwidth]{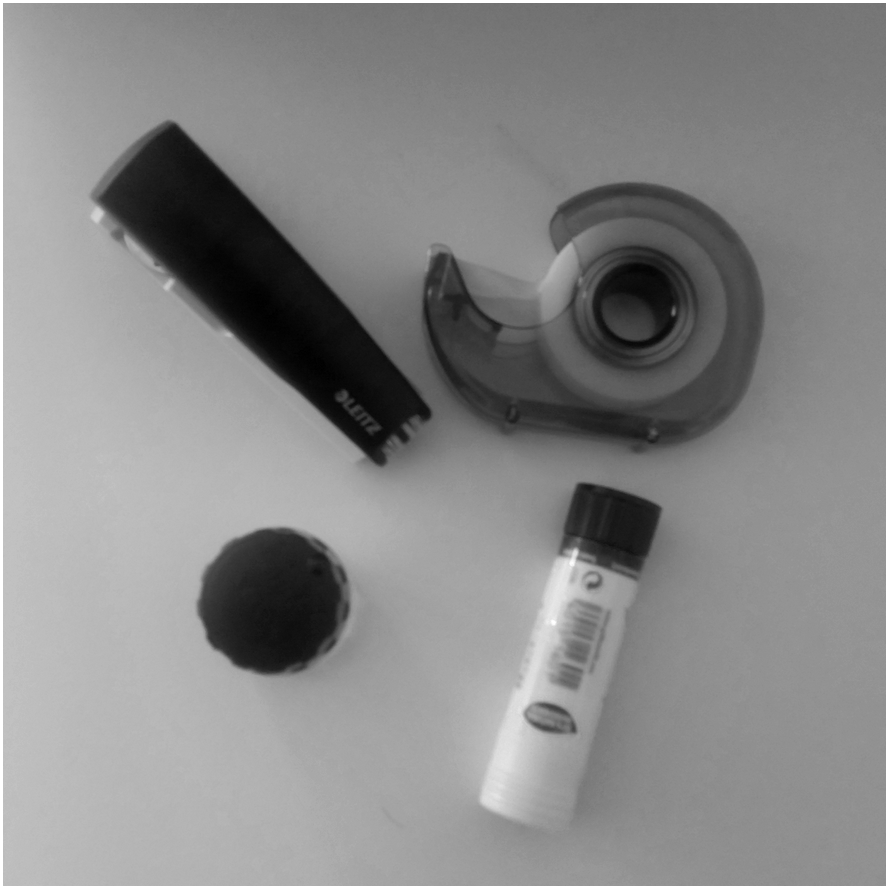}
	\includegraphics[width=0.32\textwidth]{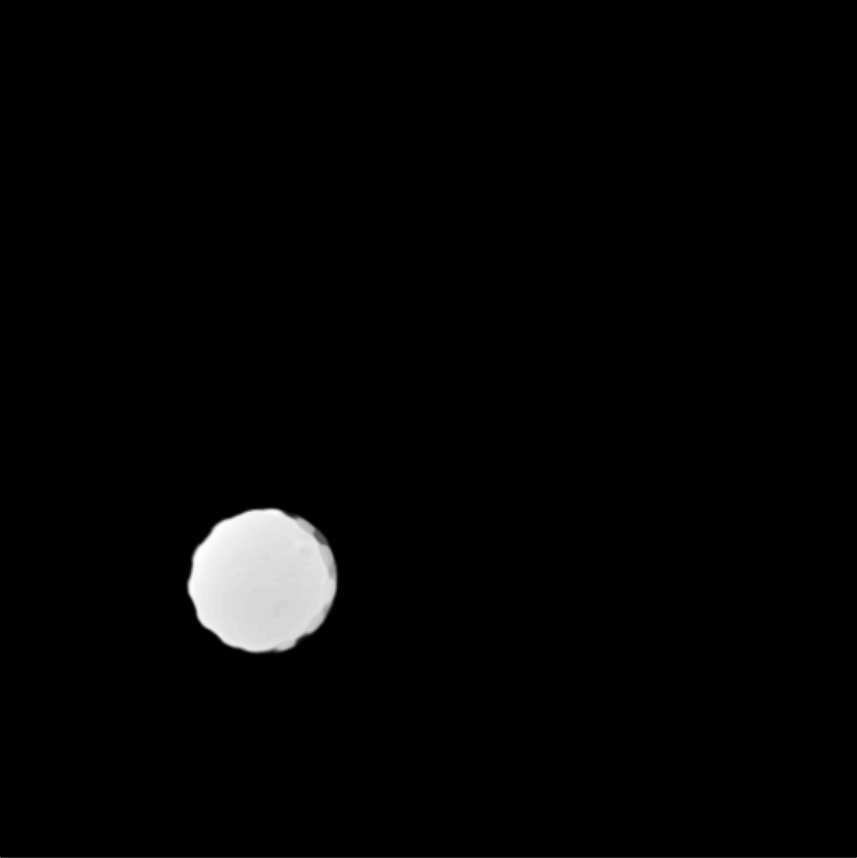}
	\includegraphics[width=0.32\textwidth]{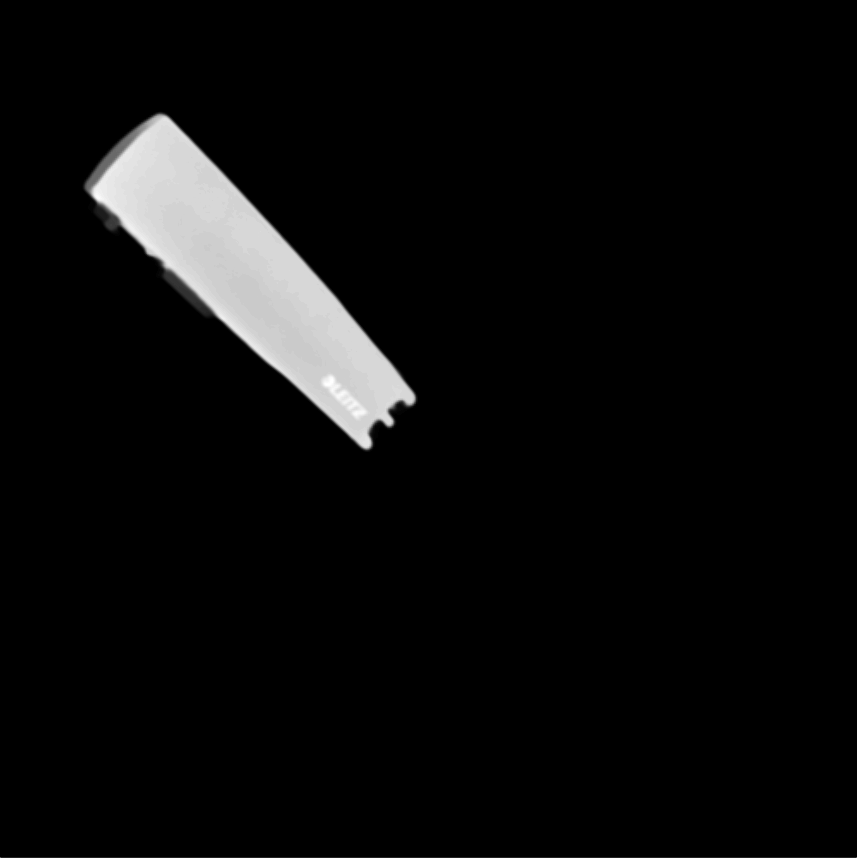}\\[0.2em]
	\includegraphics[width=0.32\textwidth]{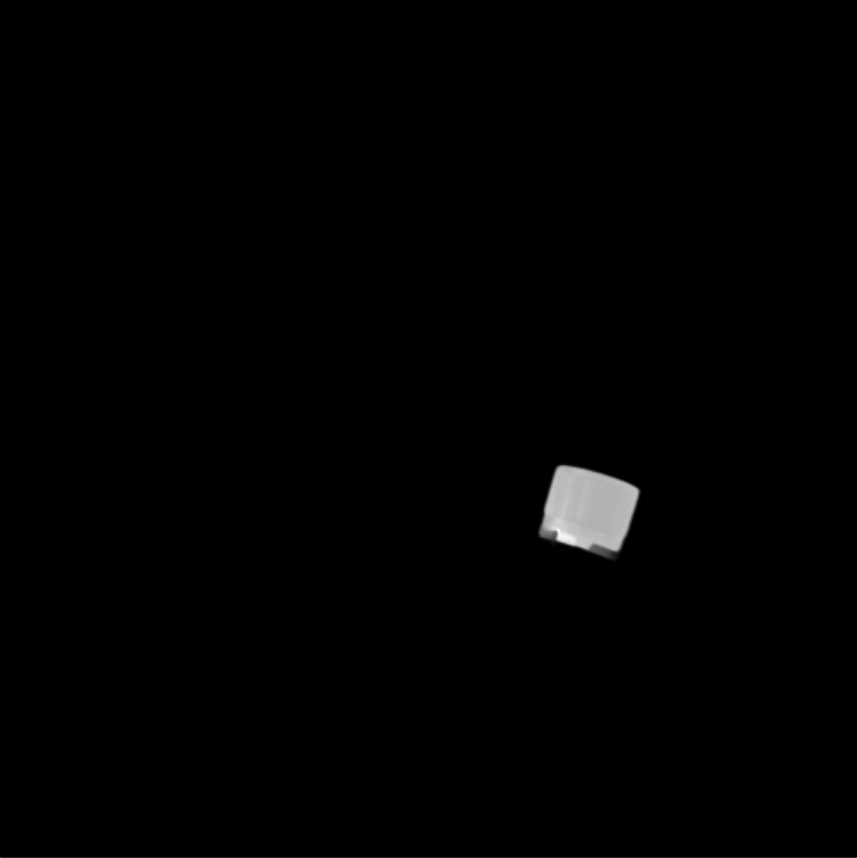}
	\includegraphics[width=0.32\textwidth]{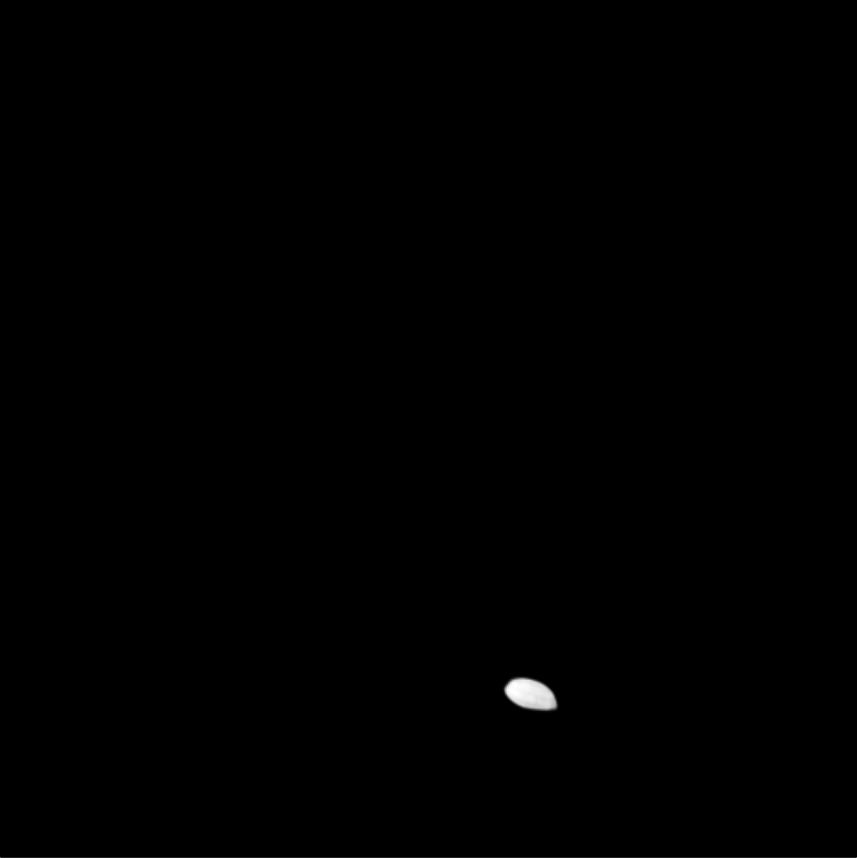}
	\includegraphics[width=0.32\textwidth]{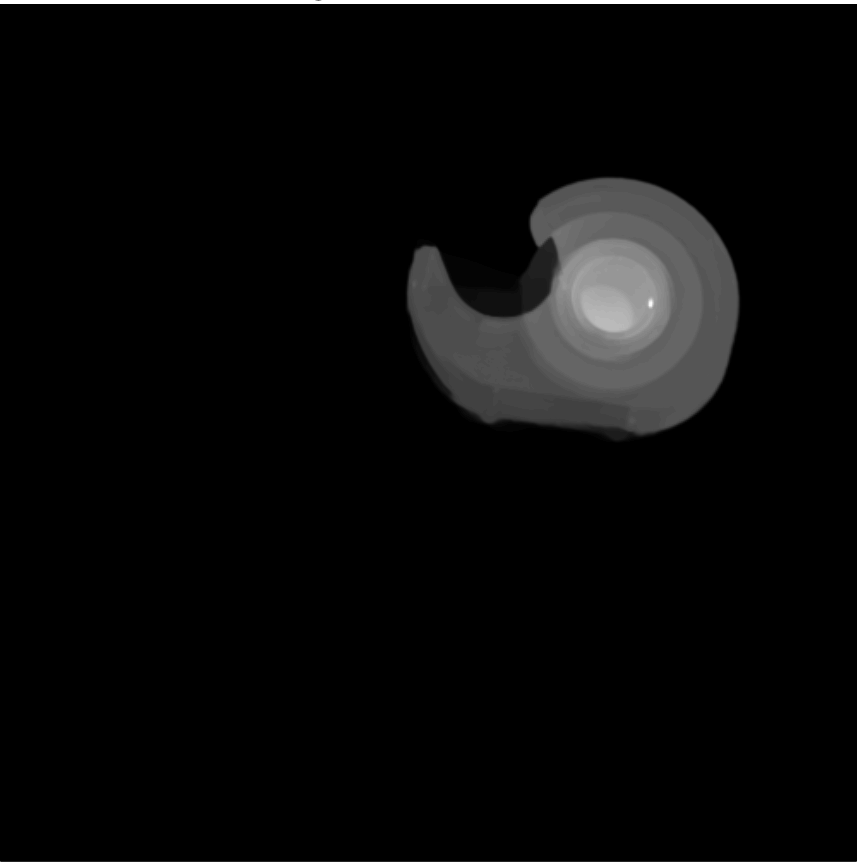}
	\caption{Two-dimensional case. From top left to bottom right: the image $I(x)$, 
	the eigenfunctions~$\phi_1,\, \phi_2,\, \phi_3,\, \phi_{4}$ and $\phi_{5}$.}
	\label{fig:things}
\end{figure}

To illustrate the remarkable approximation properties of the AE basis in two space dimensions, we now consider the 
image~$I(x)$, $x=(x_1,x_2)$, shown in Fig.~\ref{fig:things} (top, left).
Assume we want to separate the different objects appearing in the image.
We compute the first five eigenfunctions using \eqref{eq:eigenfunctionsTV} with $I(x)$ as input. Here, as in all examples in this paper, we take $\gamma=\max\left|\nabla I(x)\right|$.
These results illustrate the remarkable properties of using the AE for segmentation. We are able to extract from the image, the stapler,
the sharpener, the sellotape, the lid of the glue and the position of the slogan of the glue company.
\begin{rem}
Since the discretization of the eigenspace problem \eqref{eq:eigenfunctionsTV} is highly sparse,
and can be computed on a coarse grid, we may compute the eigenfunctions
by using a standard restarted Arnoldi iteration~\cite{Lehoucq:1996:DTI}, 
which results in a very fast algorithm ($\mathcal{O}(n)$).
The resulting eigenfunctions are highly sparse as well, see \cite{Grote:2016:AEI} for details. 
Hence, we get a good segmentation at a very low cost in terms of run-time and memory requirements.
%Moreover, the eigenfunctions are orthogonal to each other, which may be a key point for automatizing the choice of the relevant
%eigenfunctions. 
\end{rem}

	\section{Uses in Medical Imaging}
	\label{s:medical}
		
We shall now illustrate the usefulness and versatility of the method through a series of medical imaging examples.
\subsection{Segmentation of Medical Images}
	\label{ss:segmentingmedical}
\begin{figure}[ht!]
\centering	
\includegraphics[width=0.32\textwidth,height=4.2cm]{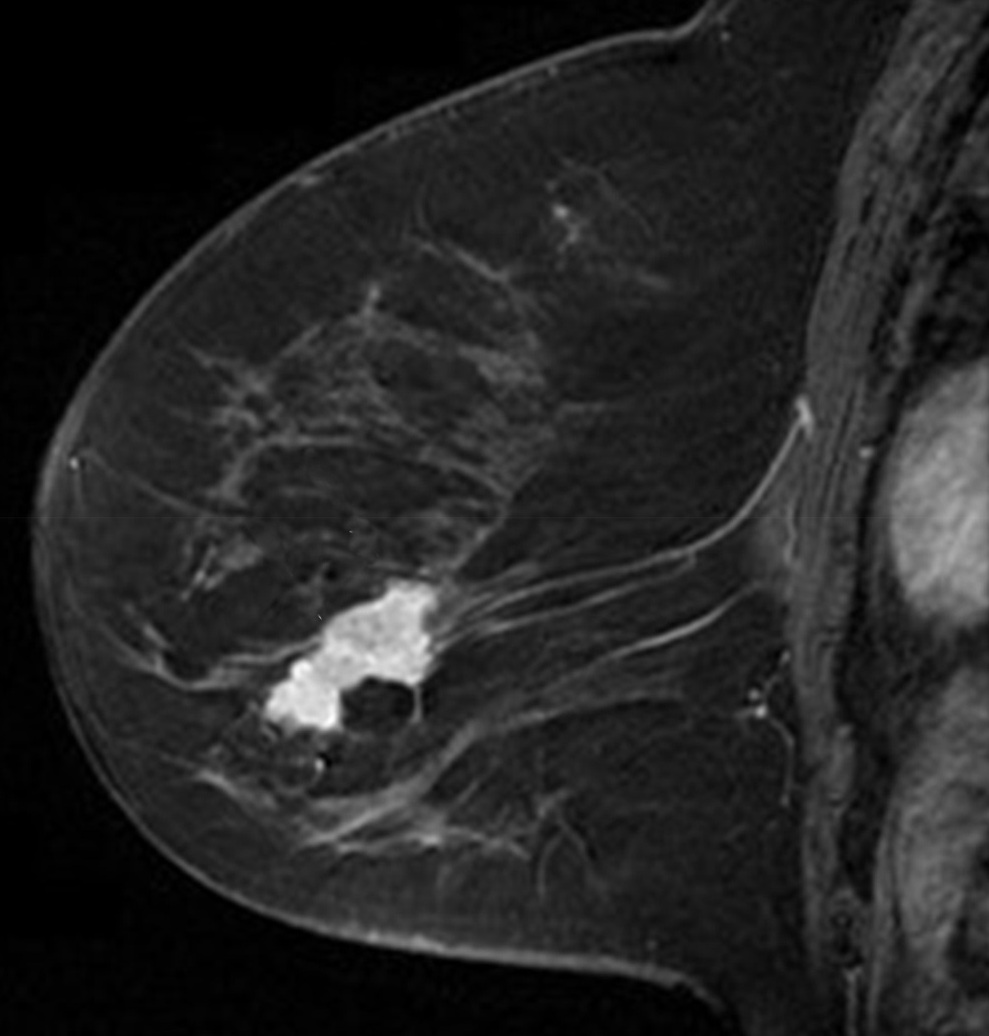}
\includegraphics[width=0.32\textwidth,height=4.2cm]{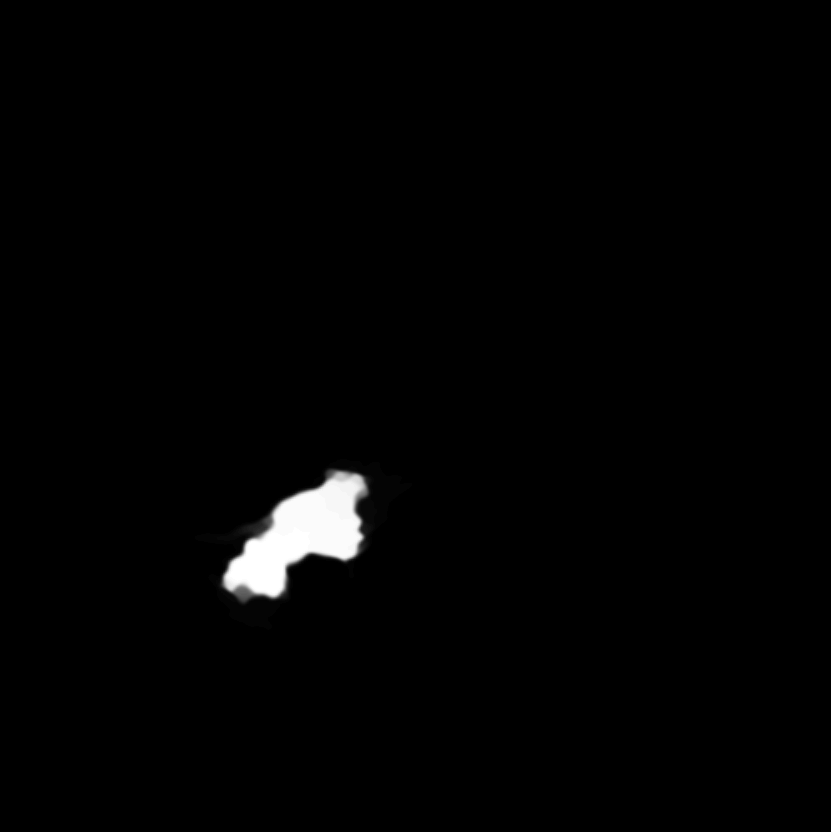}\\[0.2em]
\includegraphics[width=0.32\textwidth,height=4.2cm]{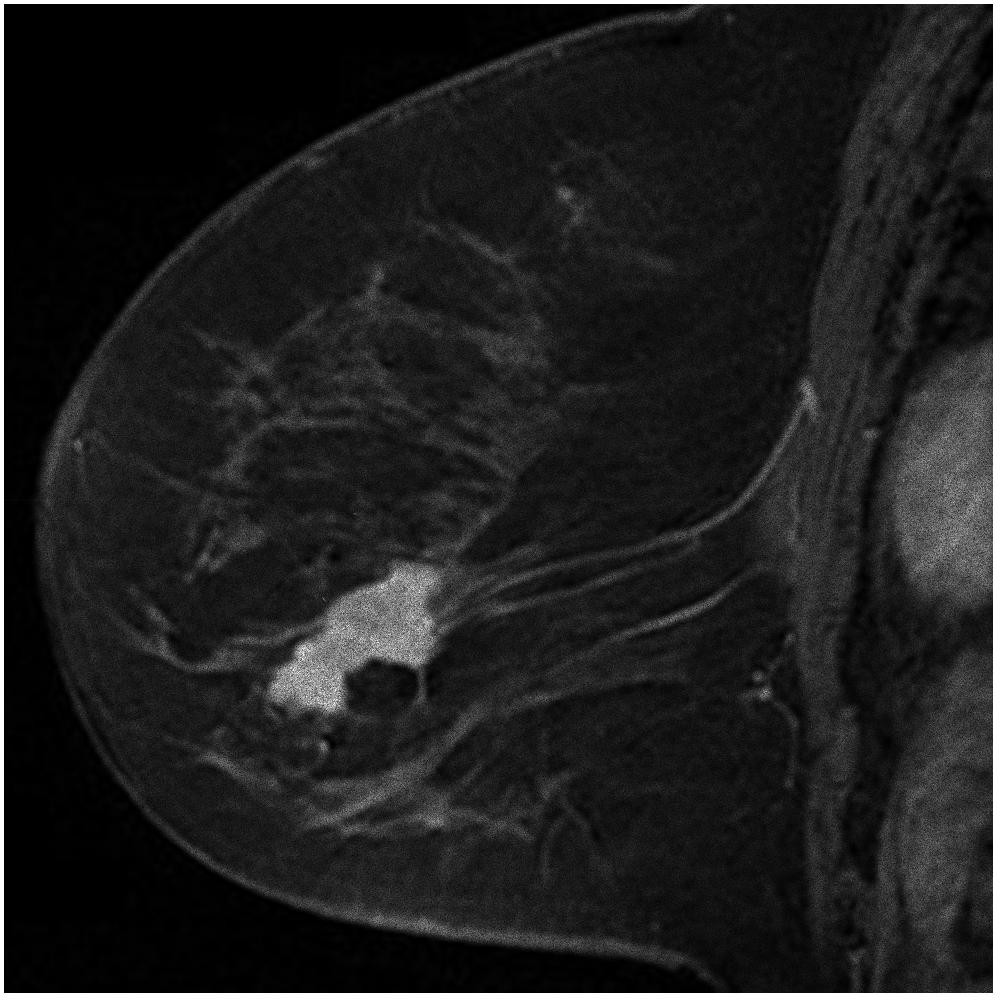}
\includegraphics[width=0.32\textwidth,height=4.2cm]{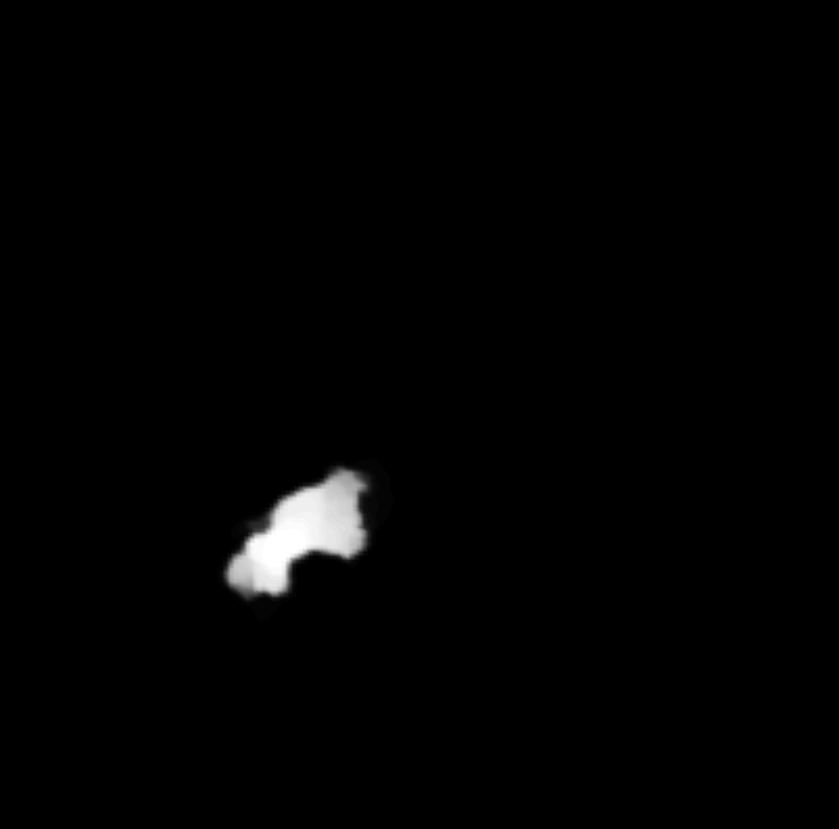}\\[0.2em]
\includegraphics[width=0.32\textwidth,height=4.2cm]{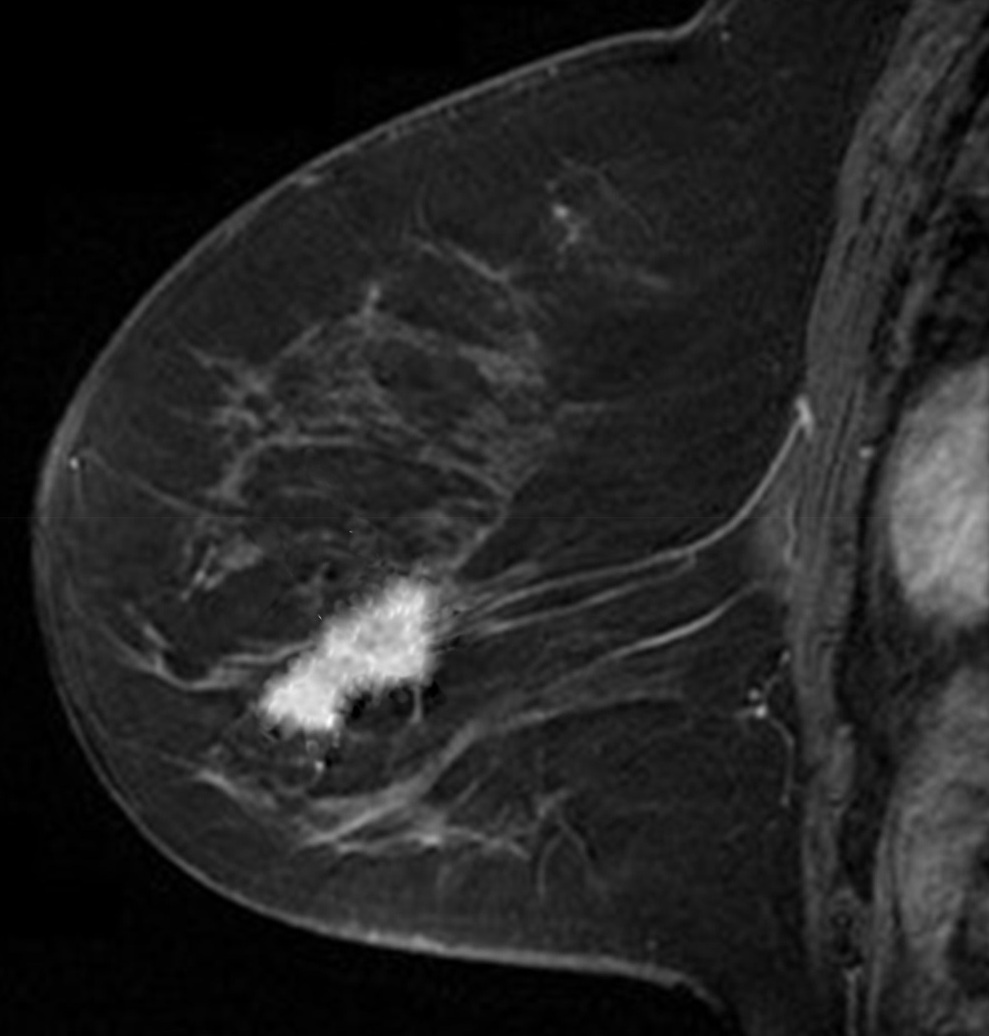}
\includegraphics[width=0.32\textwidth,height=4.2cm]{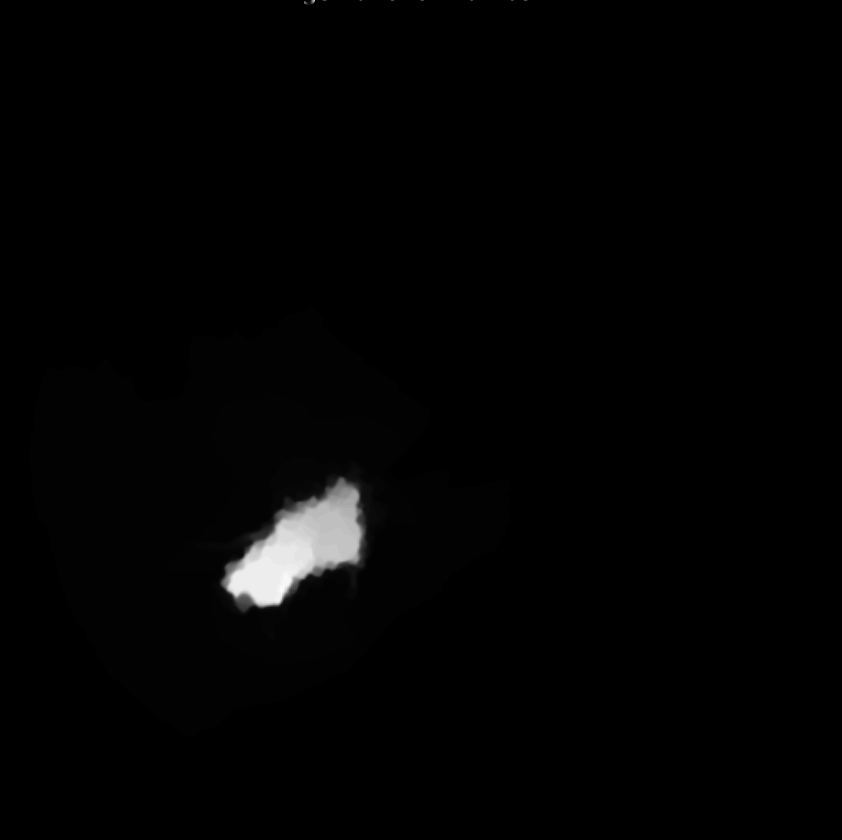}
\caption{MRI of the breast. Top: the image $I(x)$ and its first eigenfunction $\phi_1$.
Middle: the image $I_{noise}(x)$ with $20\%$ of added noise and its first eigenfunction $\phi_1$.
Bottom: the image $I_{blurred}(x)$ with blurred boundaries and its first eigenfunction $\phi_1$.}\label{fig:brust}
\end{figure}

First we will use a Magnetic Resonance (MR) image of a female breast, taken from \cite{eby2008magnetic} with permission from\footnote{\url{www.slicer.org}}, to segment the tumor mass.
On the top/right of Fig.~\ref{fig:brust}, we see the first eigenfunction segmenting the tumor perfectly out of the breast image.
Next, in the middle row of Fig.~\ref{fig:brust}, we add to the image $I$, $20\%$ of standard Gaussian noise, such that
\begin{equation}\label{noise} I_{noisy} = I\left(1+\delta\,\xi\right)\,\end{equation}
where $\xi$ is i.i.d.~Gaussian random variable with mean zero, variance equal to one and $\delta=0.2$ represents the noise level.
Now, we compute the adaptive eigenspace of $I_{noisy}$.
Again, the first eigenfunction holds a nice segmentation of the tumor despite the additional noise.
To demonstrate the quality of our approach, we consider the image $I$ but this time we destroy the boundaries of the tumor and change them by blurring (using an image manipulation program). The resulting image $I_{blurred}$ is shown on the bottom left of Fig.~\ref{fig:brust}.
On the bottom right of Fig.~\ref{fig:brust}, we see the segmentation is accurate and captures the tumor and its blurred boundaries.

\begin{figure}[h!]
	\centering	
	\includegraphics[width=0.32\textwidth,height=4.2cm]{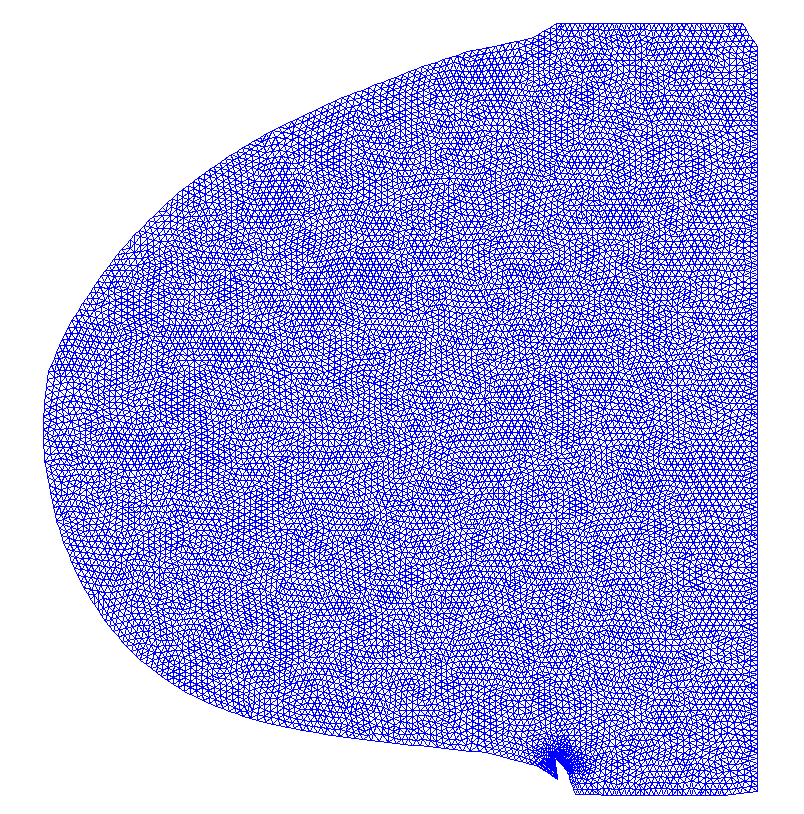}
	\includegraphics[width=0.32\textwidth,height=4.2cm]{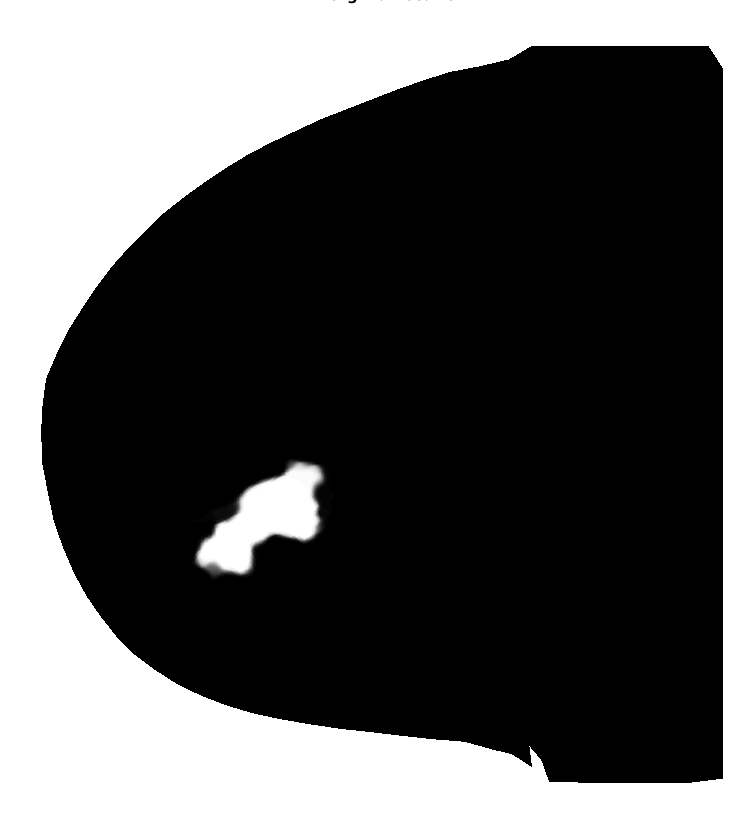}
	\includegraphics[width=0.32\textwidth,height=4.2cm]{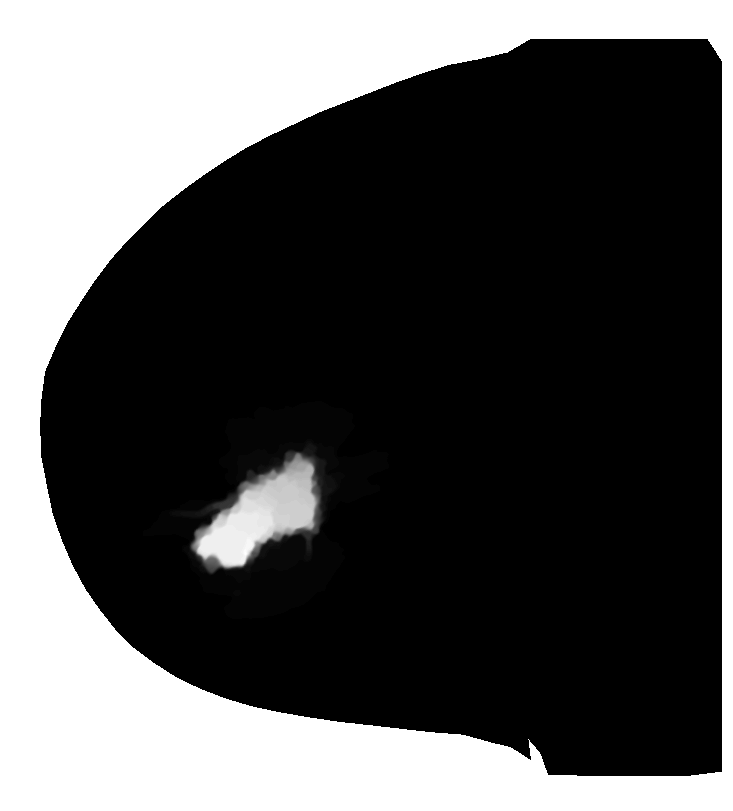}
	\caption{Segmentation of a MR image of the breast using FE: left, the automatically generated FE mesh. 
	Center: the segmentation of the tumor form image $I$. Right, the segmentation of the image $I_{blurred}$ (right) on the FE mesh.}
	\label{fig:FEbraest}
	\end{figure}

\begin{figure}[h!]
	\centering	
	\includegraphics[width=0.32\textwidth,height=4.2cm]{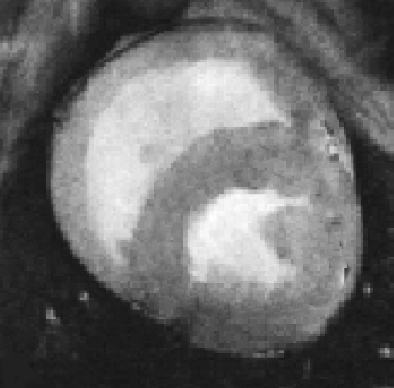}
	\includegraphics[width=0.32\textwidth,height=4.2cm]{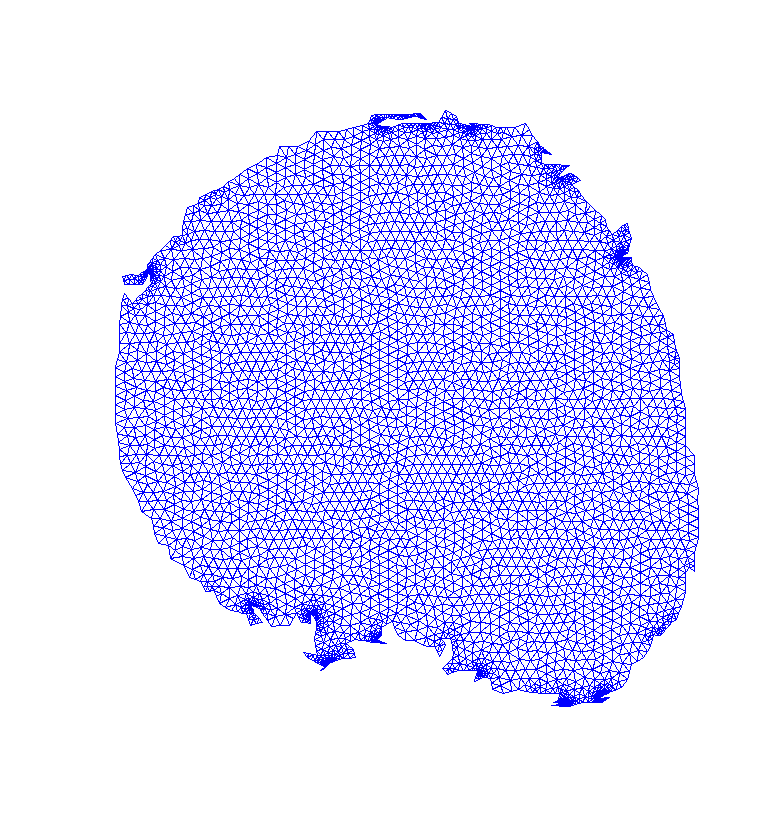}%\\[0.2em]
	\includegraphics[width=0.32\textwidth,height=4.2cm]{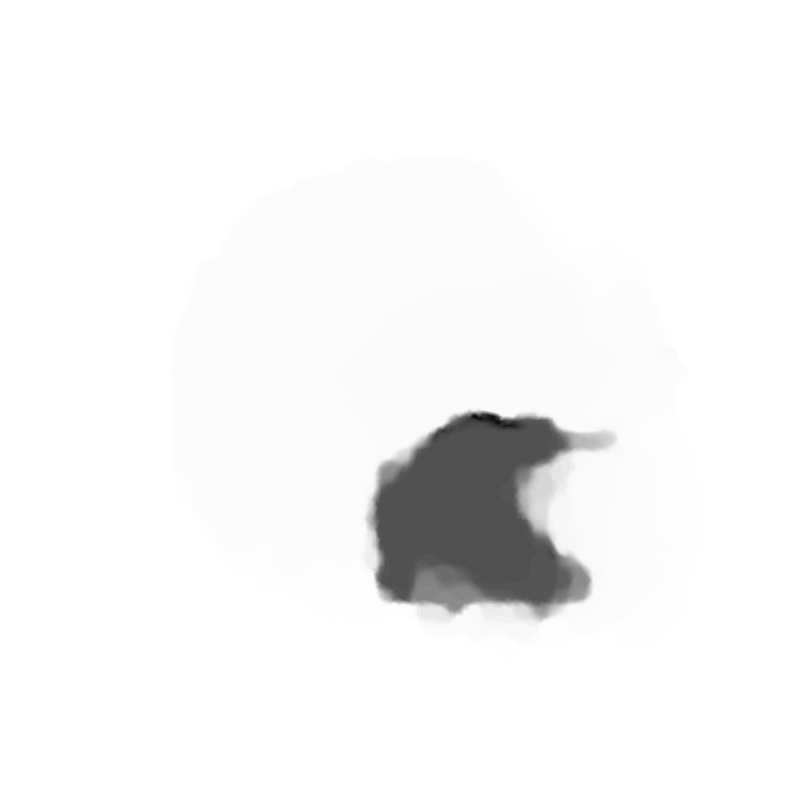}
	\caption{Segmenting MR heart image using FE: left, the MR heart image. 
	Center, the automatically generated FE mesh. Right, the first eigenfunction $\phi_1$, i.e.~the segmentation of the ventricle (right).}
	\label{fig:FEheart}
\end{figure}

To reduce computational cost and get even more accurate results, we may create a Finite Element (FE) grid on the area of interest rather than on the image boundaries. 
For illustration, we consider the MR-images $I$ and $I_{blurred}$ shown on the top/left and bottom/left of Fig.~\ref{fig:brust}, respectively. We automatically produce a mesh on the breast boundaries (see on the left of Fig.~\ref{fig:FEbraest}). The segmentation for $I$ and for $I_{blurred}$ is shown on the center and right of Fig.~\ref{fig:FEbraest}, respectively.
	This FE approach can be easily adapted to other images, for example in Fig.~\ref{fig:FEheart} left, we apply this approach for segmenting a ventricle MRI heart image, taken from \cite{Angenent2006} with permission from$^1$.
	In the center of Fig.~\ref{fig:FEheart}, the FE mesh is shown and on the right of the figure, we see the segmentation of the adaptive eigenspace using FE. As discussed in Rem.~\ref{rem:thresholding},
	we do not always get a binary segmentation, but this is easy to get using a standard threshold.
	
	\subsection{Noise Removal by Segmenting the Noise Out of an Image}
	
With the following examples we show that the adaptive eigenspace approach can not only be used to segment complicated structures but also to \textit{segment out noise and speckle} as seen in some medical images.
This can be done using \eqref{eq:param_u} with $m\leq K$, where $K$ is the number of eigenfunctions that we include in the expansion of $I$.
	Here, we take advantage of the fact that the noise/speckle in the image, appears in eigenfunctions correlating to high eigenvalues.

	Hence, we truncate the expansion of $I$ to hold only relevant information and take the first $K$ eigenfunctions related to the smallest eigenfunctions.  
	We can approximate $I$ as the following sum
	\begin{equation} \label{eq:param_u_trunc}
	\tilde{I} \ = \ I_0(x) \ + \ \sum_{m= 1}^K \beta_m \phi_m(x),
	\end{equation}
	where $I_0$, as defined in \eqref{eq:u_0}, holds the information on the boundary of the image and the eigenfunctions $\phi_m$ are computed using \eqref{eq:eigenfunctionsTV}. We have the option to set $I_0$ to zero to zero in all or part of the boundary, if we know that the boundary information there is irrelevant.

\begin{figure}[h!]
	\centering	
	\includegraphics[width=0.32\textwidth,height=4.2cm]{Octbone.png}
	\includegraphics[width=0.32\textwidth,height=4.2cm]{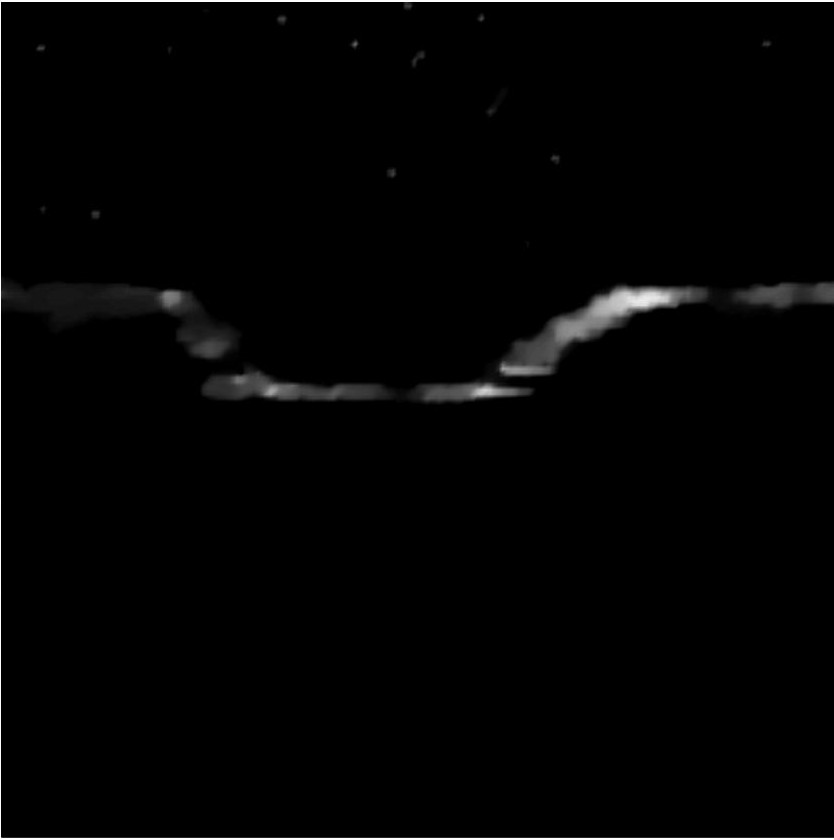}
	\caption{OCT of a bone: left, the image $I$; right The filtered image $\tilde{I}$ with~$K=150$.}
	\label{fig:octBone}
\end{figure}

We consider an OCT B-scan image of a bone piece, shown on the left of Fig.~\ref{fig:octBone}. This image is produced by measuring across the cut (yielding the cut profile)
while ablating the bone using a laser beam coming from top. In the center and right of Fig.~\ref{fig:octBone}, we see the eigenfunctions $\phi_1$ and $\phi_{250}$,
respectively. It is easy to see that $\phi_1$ is extremely relevant to the reconstruction of $\tilde{I}$ and the eigenfunction $\phi_{250}$ holding information only on
the noise. 
	\begin{figure}[ht!]
	\centering	
	\includegraphics[width=0.32\textwidth,height=4.2cm]{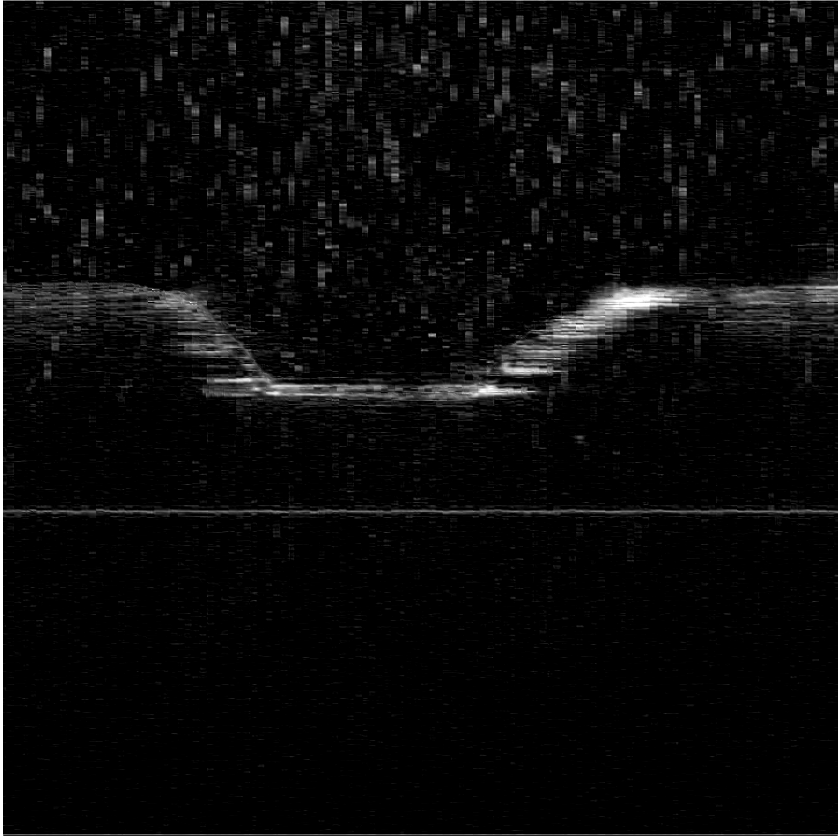} 
	\includegraphics[width=0.32\textwidth,height=4.2cm]{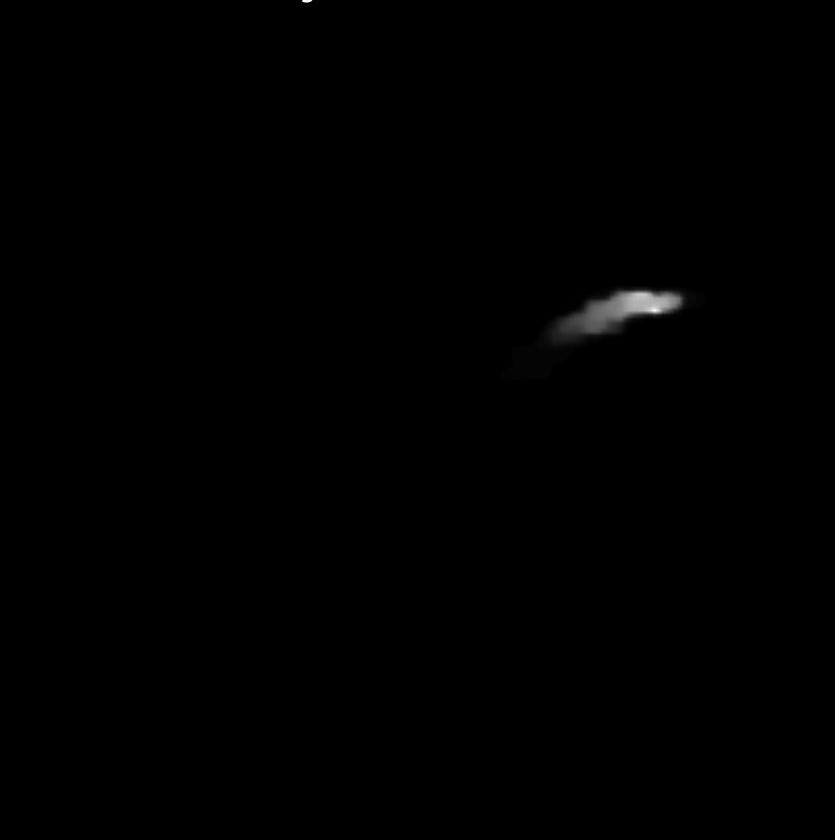}
	\includegraphics[width=0.32\textwidth,height=4.2cm]{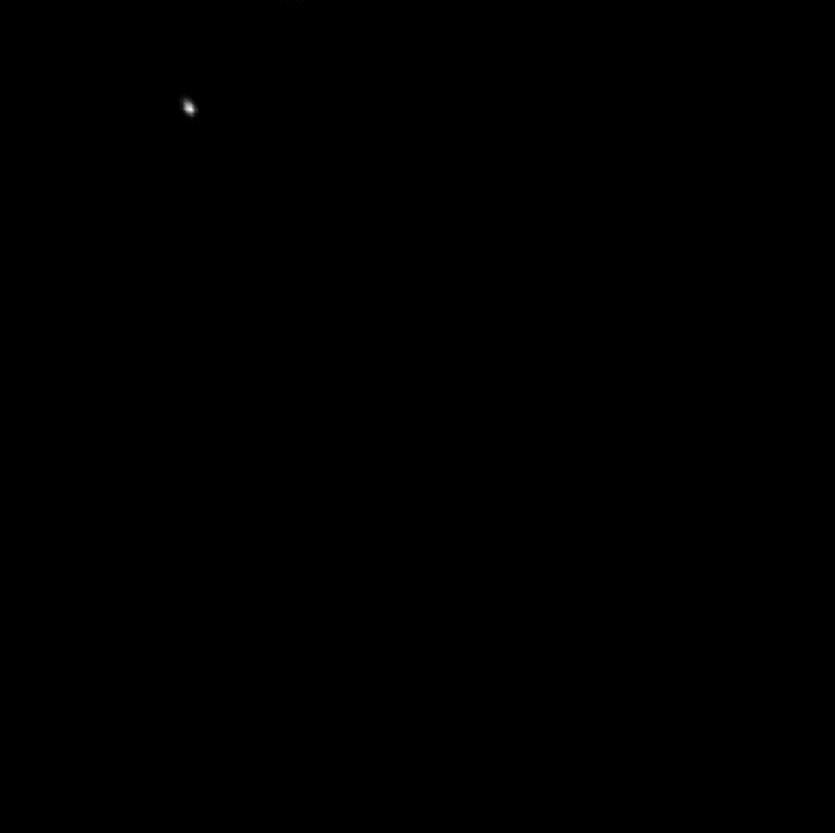}
	\caption{Image filtering: left, the image $I$. Center, the first eigenfunction~$\phi_1$. Right, the eigenfunction $\phi_{250}$.}
	\label{fig:oct}
\end{figure}

Now, we segment the noise out, produced by water droplets of the cooling spray, from the image using \eqref{eq:param_u_trunc}. The filtered image is shown on the right of Fig.~\ref{fig:octBone}.

\begin{rem}
More than $90\%$ of the entries in the eigenfunctions of the adaptive eigenspace are very small, and can be, actually, set to zero without 
loosing essential information, see~\cite{Grote:2016:AEI}. Hence, the representation of $\tilde{I}$ in \eqref{eq:param_u_trunc} is sparse and has low memory requirements. 
\end{rem}

\begin{figure}[ht!]
	\centering	
	\includegraphics[width=0.32\textwidth,height=4.2cm]{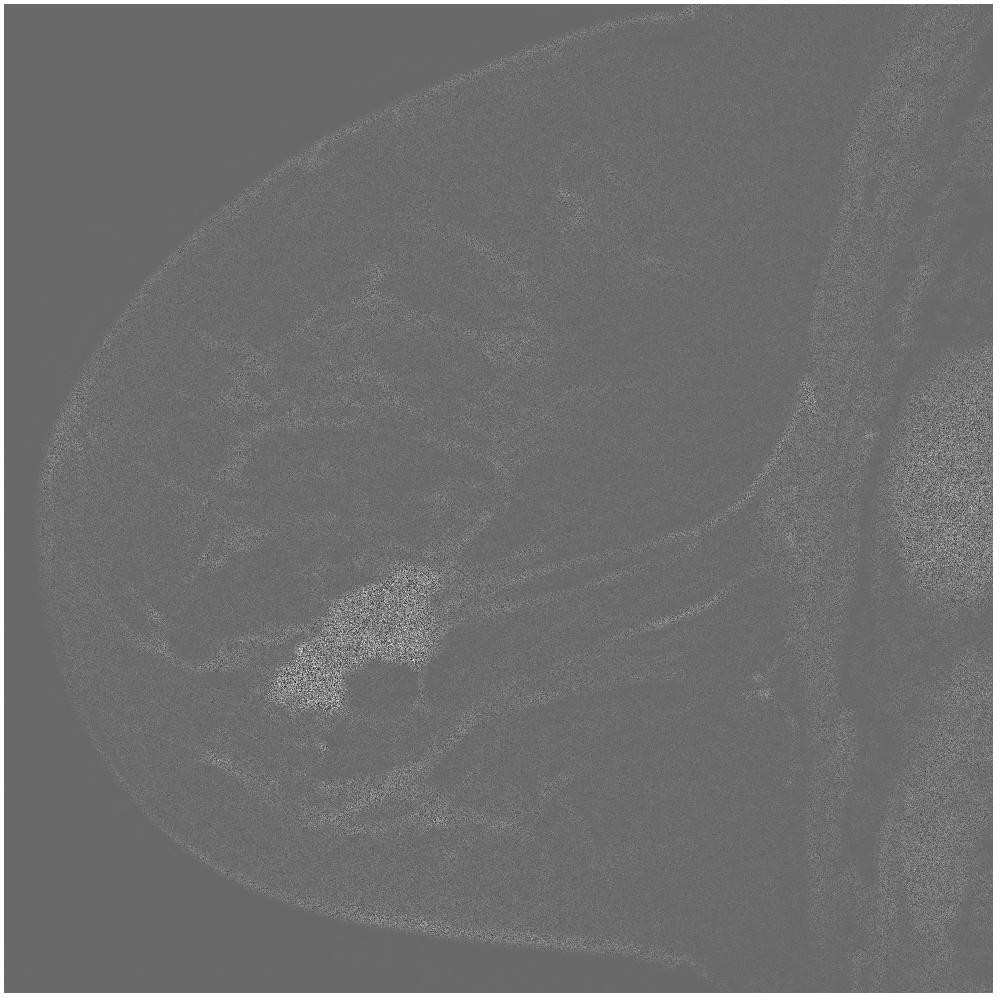}
	\includegraphics[width=0.32\textwidth,height=4.2cm]{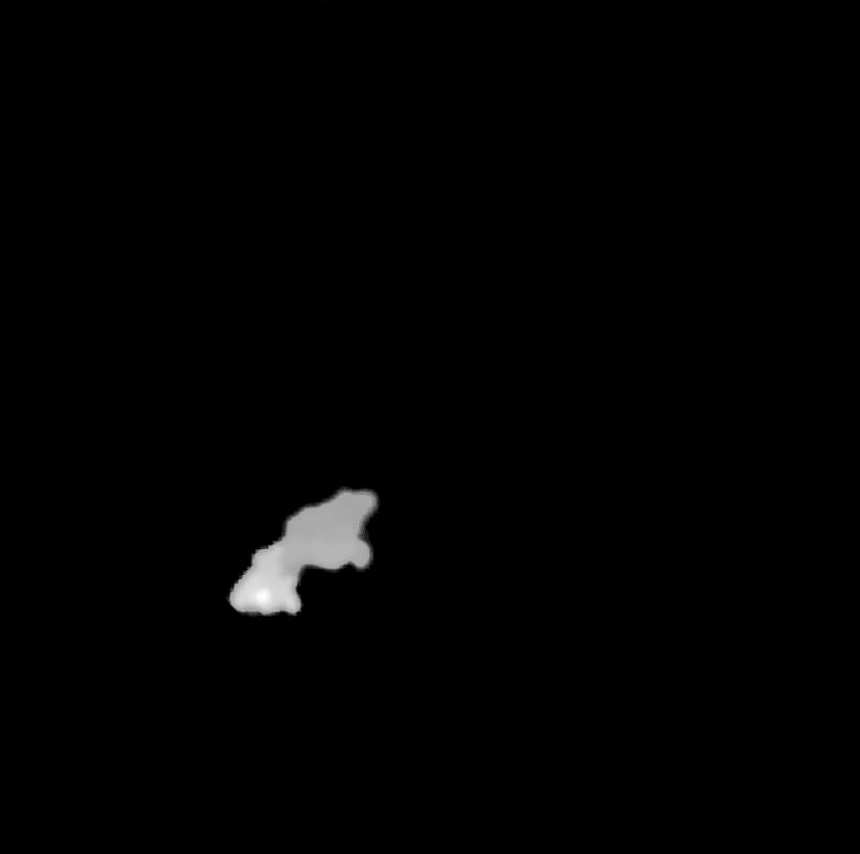}
	\caption{MRI of the breast. The image $I(x)$ with $120\%$ of added standard Gaussian noise (left) and the first eigenfunction $\phi_{1}$ of AE of the filtered image $\tilde{I}$ with $K=150$ (right).}
	\label{fig:brustnoise}
\end{figure}

In some cases the image is so noisy, that we would like to segment the noise out of it before segmenting the image. We consider once more the MR image from Fig.~\ref{fig:brust} with $120\%$ of added noise (as in \eqref{noise} with $\delta=1.2$),
such that the borders of the tumor are heavily distorted.
At the first step, we use the adaptive eigenspace to filter the image (K=150) and then we segment the tumor out of the filtered image. The segmentation is performed with much success, see on the right of Fig.~\ref{fig:brustnoise}.
The segmentation is very close to the one done on the original image, without noise, shown on Fig.~\ref{fig:brust}.

\section{Concluding Remarks}

We have presented a new framework for image segmentation based on the adaptive eigen-space. 
Instead of minimizing an energy norm and to regularize it,
we compute the eigenfunctions of the gradient of the regularization term, to segment an image. The approach has been shown to be insensitive to the parameter $\gamma=\max\left|\nabla I(x)\right|$ as the same value was used for all the experiments reported herein.
Hence,  the adaptive eigenspace segmentation is not sensitive to the parameter choice and does not need any training or other prior shape information to segment an image other than the image to segment itself.

In addition,  we showed how the adaptive eigenspace segmentation can be used to segment the noise out of an image, rather than filtering it with classical methods.

In this paper, mostly medical images are segmented, but clearly, this approach may be directly applied to other type of images.
The method uses only a fraction of the computational cost used by other segmentation methods and yet, the results are remarkable.
The eigenfunctions are highly sparse and hence this approach can be easily extended to three space dimensions.

\subsection*{Acknowledgment}
The authors thank Allen Tannenbaum for useful comments and suggestions.
%%%%%%%%%%%%%%%%%%%%%%%
%% Elsevier bibliography styles
%%%%%%%%%%%%%%%%%%%%%%%
%% To change the style, put a % in front of the second line of the current style and
%% remove the % from the second line of the style you would like to use.
%%%%%%%%%%%%%%%%%%%%%%%

%% Numbered
%\bibliographystyle{model1-num-names}

%% Numbered without titles
%\bibliographystyle{model1a-num-names}

%% Harvard
%\bibliographystyle{model2-names.bst}
%\biboptions{authoryear}

%% Vancouver numbered
%\usepackage{numcompress}
%\bibliographystyle{model3-num-names}

%% Vancouver name/year
%\usepackage{numcompress}
%\bibliographystyle{model4-names}
%\biboptions{authoryear}

%% APA style
%\bibliographystyle{model5-names}\biboptions{authoryear}

%% AMA style
%\usepackage{numcompress}
%\bibliographystyle{model6-num-names}

%% `Elsevier LaTeX' style
%\bibliographystyle{elsarticle-num}
%\bibliographystyle{abbrv}
\bibliographystyle{abbrv}
\bibliography{Refs} 

\end{document}